\documentclass[titlepage,12pt]{article}

\usepackage{amssymb}
\usepackage{latexsym}

\newcommand{\ben}{\begin{enumerate}}
\newcommand{\een}{\end{enumerate}}
\newcommand{\ble}{\begin{lem}}
\newcommand{\ele}{\end{lem}}
\newcommand{\bth}{\begin{thm}}
\renewcommand{\eth}{\end{thm}}
\newcommand{\bpr}{\begin{prop}}
\newcommand{\epr}{\end{prop}}
\newcommand{\bco}{\begin{cor}}
\newcommand{\eco}{\end{cor}}
\newcommand{\bcon}{\begin{conj}}
\newcommand{\econ}{\end{conj}}
\newcommand{\bde}{\begin{defn}}
\newcommand{\ede}{\end{defn}}
\newcommand{\bex}{\begin{exa}}
\newcommand{\eex}{\end{exa}}
\newcommand{\barr}{\begin{array}}
\newcommand{\earr}{\end{array}}
\newcommand{\btab}{\begin{tabular}}
\newcommand{\etab}{\end{tabular}}
\newcommand{\beq}{\begin{equation}}
\newcommand{\eeq}{\end{equation}}
\newcommand{\bea}{\begin{eqnarray*}}
\newcommand{\eea}{\end{eqnarray*}}
\newcommand{\bce}{\begin{center}}
\newcommand{\ece}{\end{center}}
\newcommand{\bpi}{\begin{picture}}
\newcommand{\epi}{\end{picture}}
\newcommand{\bfi}{\begin{figure} \begin{center}}
\newcommand{\efi}{\end{center} \end{figure}}
\newcommand{\capt}{\caption}
\newcommand{\bsl}{\begin{slide}{}}
\newcommand{\esl}{\end{slide}}

\newcommand{\bib}{thebibliography}
\newcommand{\pf}{{\bf Proof.}}
\newcommand{\qed}{\rule{1ex}{1ex}}
\newcommand{\Qed}{\rule{1ex}{1ex} \medskip}

\newcommand{\ol}{\overline}

\newcommand{\qmq}[1]{\quad\mbox{#1}\quad}

\newcommand{\emp}{\emptyset}
\newcommand{\sbs}{\subset}
\newcommand{\sbe}{\subseteq}

\newcommand{\spe}{\supseteq}
\newcommand{\setm}{\setminus}

\newcommand{\iso}{\cong}

\newcommand{\zh}{\hat{0}}
\newcommand{\oh}{\hat{1}}

\newcommand{\xh}{\hat{x}}
\newcommand{\lt}{\lhd}

\newcommand{\lte}{\unlhd}

\newcommand{\jn}{\vee}
\newcommand{\Jn}{\bigvee}
\newcommand{\mt}{\wedge}
\newcommand{\Mt}{\bigwedge}

\newcommand{\case}[4]{\left\{\barr{ll}#1&\mbox{#2}\\#3&\mbox{#4}\earr\right.}

\newcommand{\ce}[1]{\lceil #1 \rceil}

\newcommand{\fall}[2]{\langle{#1}\rangle_{#2}}

\newcommand{\ree}[1]{(\ref{#1})}

\newcommand{\ra}{\rightarrow}

\newcommand{\al}{\alpha}

\newcommand{\de}{\delta}

\renewcommand{\th}{\theta}
\newcommand{\ze}{\zeta}

\newcommand{\De}{\Delta}

\newcommand{\Th}{\Theta}

\newcommand{\bR}{{\bf R}}

\newcommand{\bX}{{\bf X}}
\newcommand{\bXh}{\hat{{\bf X}}}
\newcommand{\bZ}{{\bf Z}}
\newcommand{\bbC}{{\mathbb C}}
\newcommand{\bbF}{{\mathbb F}}

\newcommand{\bbR}{{\mathbb R}}
\newcommand{\bbZ}{{\mathbb Z}}
\newcommand{\cA}{{\cal A}}
\newcommand{\cB}{{\cal B}}
\newcommand{\cC}{{\cal C}}
\newcommand{\cD}{{\cal D}}
\newcommand{\cF}{{\cal F}}

\newcommand{\cR}{{\cal R}}
\newcommand{\cS}{{\cal S}}

\newcommand{\cW}{{\cal W}}

\newcommand{\cDB}{{\cal DB}}

\newcommand{\cIF}{{\cal IF}}

\newcommand{\cRF}{{\cal RF}}

\newcommand{\alt}{\tilde{\alpha}}

\newcommand{\Cb}{\ol{C}}

\newcommand{\Int}{\mathop {\rm Int}}

\newcommand{\aim}{Adv. in Math.}

\newcommand{\dm}{Discrete Math.}

\newcommand{\ejc}{European J. Combin.}
\newcommand{\im}{Invent. Math.}	

\newcommand{\jac}{J. Algebraic Combin.}

\newcommand{\jams}{J. Amer. Math. Soc.}

\newcommand{\jcta}{J. Combin. Theory Ser. A}

\newcommand{\tams}{Trans. Amer. Math. Soc.}
\newcommand{\pja}{Proc. Japan Acad. Ser. A  Math}

\newtheorem{thm}{Theorem}[section]
\newtheorem{prop}[thm]{Proposition}
\newtheorem{cor}[thm]{Corollary}
\newtheorem{lem}[thm]{Lemma}
\newtheorem{conj}[thm]{Conjecture}
\newtheorem{exa}[thm]{Example}
\newtheorem{defn}[thm]{Definition}

\newcommand{\Gaa}{\put(0,0){\circle*{3}}}

\newcommand{\Gca}{\put(20,0){\circle*{3}}}
\newcommand{\Gda}{\put(30,0){\circle*{3}}}

\newcommand{\Gga}{\put(60,0){\circle*{3}}}

\newcommand{\Gdb}{\put(30,10){\circle*{3}}}

\newcommand{\Gac}{\put(0,20){\circle*{3}}}

\newcommand{\Gdc}{\put(30,20){\circle*{3}}}
\newcommand{\Gec}{\put(40,20){\circle*{3}}}

\newcommand{\Ggc}{\put(60,20){\circle*{3}}}

\newcommand{\Gad}{\put(0,30){\circle*{3}}}

\newcommand{\Gdd}{\put(30,30){\circle*{3}}}

\newcommand{\Ggd}{\put(60,30){\circle*{3}}}

\newcommand{\Gae}{\put(0,40){\circle*{3}}}

\newcommand{\Gce}{\put(20,40){\circle*{3}}}
\newcommand{\Gde}{\put(30,40){\circle*{3}}}

\newcommand{\Gge}{\put(60,40){\circle*{3}}}

\newcommand{\Gdf}{\put(30,50){\circle*{3}}}

\newcommand{\Gag}{\put(0,60){\circle*{3}}}

\newcommand{\Gdg}{\put(30,60){\circle*{3}}}
\newcommand{\Geg}{\put(40,60){\circle*{3}}}

\newcommand{\GaaL}[2]{\Gaa \put(0,0){\makebox(0,0){\hspace{#1}#2}}}

\newcommand{\GcaL}[2]{\Gca \put(20,0){\makebox(0,0){\hspace{#1}#2}}}

\newcommand{\GgaL}[2]{\Gga \put(60,0){\makebox(0,0){\hspace{#1}#2}}}

\newcommand{\GacL}[2]{\Gac \put(0,20){\makebox(0,0){\hspace{#1}#2}}}

\newcommand{\GdcL}[2]{\Gdc \put(30,20){\makebox(0,0){\hspace{#1}#2}}}
\newcommand{\GecL}[2]{\Gec \put(40,20){\makebox(0,0){\hspace{#1}#2}}}

\newcommand{\GgcL}[2]{\Ggc \put(60,20){\makebox(0,0){\hspace{#1}#2}}}

\newcommand{\GadL}[2]{\Gad \put(0,30){\makebox(0,0){\hspace{#1}#2}}}

\newcommand{\GddL}[2]{\Gdd \put(30,30){\makebox(0,0){\hspace{#1}#2}}}

\newcommand{\GgdL}[2]{\Ggd \put(60,30){\makebox(0,0){\hspace{#1}#2}}}

\newcommand{\GaeL}[2]{\Gae \put(0,40){\makebox(0,0){\hspace{#1}#2}}}

\newcommand{\GceL}[2]{\Gce \put(20,40){\makebox(0,0){\hspace{#1}#2}}}
\newcommand{\GdeL}[2]{\Gde \put(30,40){\makebox(0,0){\hspace{#1}#2}}}

\newcommand{\GgeL}[2]{\Gge \put(60,40){\makebox(0,0){\hspace{#1}#2}}}

\newcommand{\GagL}[2]{\Gag \put(0,60){\makebox(0,0){\hspace{#1}#2}}}

\newcommand{\GdgL}[2]{\Gdg \put(30,60){\makebox(0,0){\hspace{#1}#2}}}
\newcommand{\GegL}[2]{\Geg \put(40,60){\makebox(0,0){\hspace{#1}#2}}}

\newcommand{\Gaaga}{\put(0,0){\line(1,0){60}}}

\newcommand{\Gaaac}{\put(0,0){\line(0,1){20}}}

\newcommand{\Gaadg}{\put(0,0){\line(1,2){30}}}

\newcommand{\Gcaac}{\put(20,0){\line(-1,1){20}}}

\newcommand{\Gcaec}{\put(20,0){\line(1,1){20}}}

\newcommand{\Gdaac}{\put(30,0){\line(-3,2){30}}}

\newcommand{\Gdadc}{\put(30,0){\line(0,1){20}}}

\newcommand{\Gdagc}{\put(30,0){\line(3,2){30}}}

\newcommand{\Ggadg}{\put(60,0){\line(-1,2){30}}}

\newcommand{\Gdbad}{\put(30,10){\line(-3,2){30}}}

\newcommand{\Gdbdd}{\put(30,10){\line(0,1){20}}}

\newcommand{\Gdbgd}{\put(30,10){\line(3,2){30}}}

\newcommand{\Gacae}{\put(0,20){\line(0,1){20}}}

\newcommand{\Gacce}{\put(0,20){\line(1,1){20}}}
\newcommand{\Gacde}{\put(0,20){\line(3,2){30}}}

\newcommand{\Gdcae}{\put(30,20){\line(-3,2){30}}}

\newcommand{\Gdcge}{\put(30,20){\line(3,2){30}}}

\newcommand{\Gecce}{\put(40,20){\line(-1,1){20}}}

\newcommand{\Gecge}{\put(40,20){\line(1,1){20}}}

\newcommand{\Ggcde}{\put(60,20){\line(-3,2){30}}}

\newcommand{\Ggcge}{\put(60,20){\line(0,1){20}}}

\newcommand{\Gaddf}{\put(0,30){\line(3,2){30}}}

\newcommand{\Gdddf}{\put(30,30){\line(0,1){20}}}

\newcommand{\Ggddf}{\put(60,30){\line(-3,2){30}}}

\newcommand{\Gaeag}{\put(0,40){\line(0,1){20}}}

\newcommand{\Gaedg}{\put(0,40){\line(3,2){30}}}

\newcommand{\Gceeg}{\put(20,40){\line(1,1){20}}}

\newcommand{\Gdedg}{\put(30,40){\line(0,1){20}}}

\newcommand{\Ggedg}{\put(60,40){\line(-3,2){30}}}
\newcommand{\Ggeeg}{\put(60,40){\line(-1,1){20}}}

\begin{document}
\pagestyle{empty}
\title{Why the characteristic polynomial factors
}
\author{Bruce E. Sagan\\
Department of Mathematics\\ 
Michigan State University\\
East Lansing, MI 48824-1027\\
sagan@math.msu.edu}


\date{\today \\[1in]
	\begin{flushleft}
	\end{flushleft}
       }
\maketitle

\begin{flushleft} Proposed running head: \end{flushleft}
	\begin{center} 
Why $\chi$ factors
	\end{center}

Send proofs to:
\begin{center}
Bruce E. Sagan \\ Department of Mathematics \\Michigan State
University \\ East Lansing, MI 48824-1027\\[5pt]
Tel.: 517-355-8329
FAX: 517-336-1562\\
Email: sagan@mth.msu.edu
\end{center}

	\begin{abstract}
We survey three methods for proving that the characteristic polynomial
of a finite ranked lattice factors over the nonnegative integers and
indicate how they have evolved recently.  The first 
technique uses geometric ideas and 
is based on Zaslavsky's theory of signed graphs.
The second approach is algebraic and employs results of Saito and
Terao about free hyperplane arrangements.  Finally we consider 
a purely combinatorial theorem of Stanley about supersolvable lattices
and its generalizations.
	\end{abstract}
\pagestyle{plain}

\section{Introduction}					\label{i}

The fundamental problem in enumerative combinatorics can be stated:
given a set $S$, find a formula for its cardinality $|S|$.  More
generally, given a sequence of sets $S_0,S_1,S_2,\ldots$ we would
like to investigate properties of the sequence 
\beq							\label{a}
a_0,a_1,a_2,\ldots
\eeq
where $a_i=|S_i|,\ i\ge0$.  From their definition we obviously have
$a_i\in\bbZ_{\ge0}$, the non-negative integers.

We can also turn this problem around.  Suppose we are given a
sequence~(\ref{a}) where the $a_i$ are defined in a way that would
permit them to be complex numbers, $\bbC$.  If, in fact, the
terms are in 
$\bbZ_{\ge0}$, can we find a combinatorial explanation for this?
One possibility, of course, would be to find a sequence of sets such
that $a_i=|S_i|$.  Such questions are of great interest currently in
algebraic combinatorics. 

We are going to investigate a particular problem of this type: trying
to explain why the roots of a certain polynomial associated with a
partially ordered set, called the characteristic polynomial, often has all
of its roots in $\bbZ_{\ge0}$.   We will provide three explanations with tools
drawn from three different areas of mathematics: graph
theory/geometry, algebra, and pure combinatorics.   The first of these
uses Zaslavsky's lovely theory of signed
graph coloring~\cite{zas:grs,zas:sgc,zas:cis} which can be generalized
to counting points of $\bbZ^n$ or of $\bbF_p^n$ inside a certain 
polytope~\cite{ath:cps,bs:cep,cr:cg,hai:cqr,ter:jdf}.
(Here $\bbF_p$ is the Galois field with $p$ elements.)
The next technique is based on theorems of Saito~\cite{sai:tld} and
Terao~\cite{ter:gef} about free hyperplane arrangements.  Work has
also been done on related concepts such as inductive
freeness~\cite{ter:ahf} and recursive freeness~\cite{zie:ach}.
The third method employs a theorem of Stanley~\cite{sta:sl} on
semimodular supersolvable lattices which has recently been generalized
by Blass and myself~\cite{bs:mfl} by relaxing both restrictions on the lattice.
Along the way we will meet many important combinatorial concepts.

The rest of this paper is organized as follows.  The next section
will introduce the M\"obius function, $\mu$, of a partially ordered
set (poset) which is a far-reaching generalization of the one in number theory.
In section~\ref{gfcp} we will talk about generating functions, an
important way to manipulate sequences such as~(\ref{a}), and define
the characteristic polynomial, $\chi$, as the generating function for $\mu$.
The last three sections will be devoted to the three methods for
proving that for various posets $\chi$ factors over $\bbZ_{\ge0}$.

\section{M\"obius functions and posets}

In number theory, one usually sees the M\"obius
function $\mu:\bbZ_{>0}\ra\bbZ$ defined as
\beq						\label{nt}
\mu(n)=\case{0}{if $n$ is not square free,}{(-1)^k}{if $n$ is the
product of $k$ distinct primes,}
\eeq
a definition which seems very strange at first blush.  The importance
of $\mu$ lies in the number-theoretic M\"obius Inversion Theorem.
\bth						\label{ntmit}
Let $f,g:\bbZ_{>0}\ra\bbZ$ satisfy  
$$
f(n)=\sum_{d|n}g(d)
$$
for all $n\in\bbZ_{>0}$.  Then
$$
g(n)=\sum_{d|n}\mu(n/d)f(d).\qquad\qed
$$
\eth

Moving into the area of enumerative combinatorics, one of the 
very useful tools is the
Principle of Inclusion-Exclusion or PIE.
\bth						\label{pie}
Let $S$ be a finite set and $S_1,\ldots,S_n\sbe S$ then
$$
|S-\bigcup_{i=1}^n S_i|=|S|-\sum_{1\le i\le n} |S_i|+
	\sum_{1\le i<j\le n} |S_i\cap S_j|-\cdots+
	(-1)^n|\bigcap_{i=1}^n S_i|.\qquad\qed
$$
\eth

In the theory of difference equations if one takes a function 
$f:\bbZ_{\ge0}\ra\bbC$ then there is an analog of the derivative,
namely the difference operator 
$$
\De f(n)=f(n)-f(n-1)
$$ 
(where $f(-1)=0$ by definition), and an analog
of the integral, namely the summation operator 
$$
S f(n)=\sum_{i=0}^n f(i).
$$
The Fundamental Theorem of the Difference
Calculus then states
\bth						\label{ftdc}
If $f:\bbZ_{\ge0}\ra\bbC$ then
$$
\De S f(n)=f(n).\qquad\qed
$$
\eth

One of the advantages of the combinatorial M\"obius function is that its
inversion theorem unifies and generalizes the previous three results.
In addition, it makes the definition~(\ref{nt}) transparent, encodes 
topological information about posets~\cite{bjo:hsm,rot:tmf}, and has
even been used to bound the running time of certain
algorithms~\cite{bly:ldt}.  We will now define this powerful invariant.

Let finite $P$ be a poset with partial order $\le$. If $P$ has a unique
minimal element then it will be denoted $\zh=\zh_P$, and if it has a unique
maximal element then we will use the notation $\oh=\oh_P$.
If $x\le y$ in $P$ then the corresponding {\it (closed) interval} is
$$
[x,y]=\{z\ :\ x\le z\le y\}
$$
and we let $\Int(P)$ denote the set of all intervals of $P$,
Note that $[x,y]$ is a poset in its own right with $\zh_{[x,y]}=x,
\oh_{[x,y]}=y$.  
The {\it M\"obius function of $P$}, $\mu: \Int(P)\ra \bZ$, is defined
recursively by
\beq					\label{mu}
\mu(x,y)=\case{1}{if $x=y$,}{-\sum_{x\le z<y}\mu(x,z)}{if $x<y$.}
\eeq
Equivalently
\beq					\label{mu2}
\sum_{x\le z\le y} \mu(x,z)=\de_{x,y}
\eeq
where $\de_{x,y}$ is the Kronecker delta.  If $P$ has a zero
then we define $\mu(x)=\mu(\zh,x)$.

\thicklines
\setlength{\unitlength}{2pt}
\bfi
\btab{lll}
\bpi(20,90)(-10,-20)
\GaaL{-20pt}{0}
\GacL{-20pt}{1}
\GaeL{-20pt}{2}
\GagL{-20pt}{3}
\Gaaac \Gacae \Gaeag
\put(0,-20){\makebox(0,0){The chain $C_3$}}
\epi&
\bpi(80,90)(-10,-20)
\Gda
\put(30,-10){\makebox(0,0){$\emp$}}
\put(30,-20){\makebox(0,0){The Boolean algebra $B_3$}}
\GacL{-25pt}{$\{1\}$}
\GdcL{-25pt}{$\{2\}$}
\GgcL{-25pt}{$\{3\}$}
\GaeL{-35pt}{$\{1,2\}$}
\GdeL{-35pt}{$\{1,3\}$}
\GgeL{-35pt}{$\{2,3\}$}
\Gdg
\put(30,70){\makebox(0,0){$\{1,2,3\}$}}
\Gdaac \Gdadc \Gdagc \Gacae \Gacde \Gdcae \Gdcge \Ggcde \Ggcge \Gaedg 
\Gdedg \Ggedg
\epi&
\bpi(80,90)(-10,-20)
\GcaL{20pt}{1}
\GacL{-20pt}{2}
\GecL{20pt}{3}
\GceL{-20pt}{6}
\GgeL{20pt}{9}
\GegL{-20pt}{18}
\Gcaac \Gcaec \Gacce \Gecce \Gecge \Gceeg \Ggeeg
\put(30,-20){\makebox(0,0){The divisor poset $D_{18}$}}
\epi
\etab
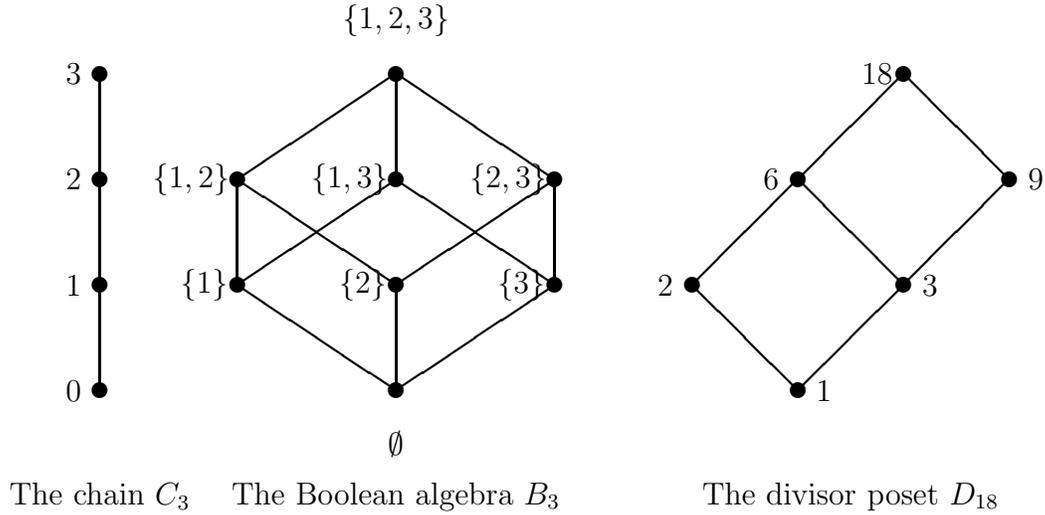
\capt{Some example posets}  \label{ex}
\efi

Let us compute $\mu(x)$ in some simple posets.
The {\it chain}, $C_n$,  consists of
the integers $\{0,1,\ldots,n\}$ ordered in the usual manner; see
Figure~\ref{ex} for a picture of $C_3$.  It is immediate directly from
the definition~(\ref{mu}) that in $C_n$ we have 
\beq					\label{muC}
\mu(x)={\left\{\barr{ll}1&\mbox{if $x=0$}\\
	 -1&\mbox{if $x=1$}\\
	  0&\mbox{if $x\ge2$}
	\earr\right.}
\eeq

The {\it Boolean algebra}, $B_n$, has as elements all subsets of
$[n]:=\{1,2,\ldots,n\}$ and $\sbe$ as order relation, the case $n=3$ being
displayed in Figure~\ref{ex}.  After computing the M\"obius function
of $B_3$, the reader will immediately guess that if $x\in B_n$ then
$\mu(x)=(-1)^{|x|}$.  This follows easily from the following
observations.  The Cartesian product of two posets $P,Q$ is
obtained by ordering the $(x,y)\in P\times Q$ component-wise.  It is
easy to prove directly from~(\ref{mu}) that if $\zh_P$ 
and $\zh_Q$ exist then
$$
\mu_{P\times Q}(x,y)=\mu_P(x)\mu_Q(y)
$$
Since $B_n$ is isomorphic as a poset to the $n$-fold product
$(C_1)^n$, it is a simple matter to verify that its M\"obius function
has the desired form using equation~\ree{muC}.

The {\it divisor poset}, $D_n$, consists of all $d|n$ ordered by
$c\le d$ if $c|d$.  Figure~\ref{ex} shows $D_{18}$.  Clearly if
$n$ has prime factorization $n=\prod_i p_i^{n_i}$ then
we have an isomorphism $D_n\iso \times_{i} C_{n_i}$.  So as with the
Boolean algebra, we get the value of $\mu(d)$ to be as in
definition~(\ref{nt}), this time in a much more natural way.

In case the reader is not convinced that the definition~\ree{mu2} is
natural, consider the {\it incidence algebra of $P$}, $I(P)$, which
consists of all functions $f:\Int(P)\ra\bbC$.  The multiplication in
this algebra is {\it convolution} defined by
$$
f*g(x,y)=\sum_{x\le z\le y} f(x,z)g(z,y).
$$
Note that with this multiplication $I(P)$ has an identity element
$\de:\Int(P)\ra\bbC$, namely $\de(x,y)=\de_{x,y}$.
One of the simplest but most important functions in $I(P)$ is the 
{\it zeta function} (so called because in $I(D_n)$ it is related to
the Riemann zeta function) given by $\ze(x,y)=1$ for all intervals
$[x,y]$.  It is easy to see that $\ze$ is invertible in $I(P)$ and in
fact that $\ze^{-1}=\mu$ where $\mu$ is defined by~\ree{mu2}.

The fundamental result about $\mu$ is the combinatorial M\"obius
Inversion Theorem~\cite{rot:tmf}.
\bth						\label{mit}
Let $P$ be a finite poset and $f,g:P\ra\bbC$.
\ben
\item  If for all $x\in P$ we have $f(x)=\sum_{y\le x} g(y)$ then
$$
g(x)=\sum_{y\le x}\mu(y,x)f(y).
$$
\item  If for all $x\in P$ we have $f(x)=\sum_{y\ge x} g(y)$ then
$$
g(x)=\sum_{y\ge x}\mu(x,y)f(y). \qquad\qed
$$
\een
\eth
It is now easy to obtain the Theorems~\ref{ntmit}, \ref{pie}, and \ref{ftdc}
as corollaries by using M\"obius inversion over $D_n$, $B_n$, and
$C_n$, respectively.  For example, to get the Principle of
Inclusion-Exclusion, use $f,g:B_n\ra\bbZ_{\ge0}$ defined by
\bea
f(X)&=&|S_X|,\\
g(X)&=&|S_X-\bigcup_{i\not\in X}S_i|,
\eea
where $S_X=\cap_{i\in X} S_i$.

\section{Generating functions and characteristic polynomials} \label{gfcp}

The {\it (ordinary) generating function} for the sequence~(\ref{a}) is
the formal power series
$$
f(x)=a_0+a_1x+a_2x^2+\cdots.
$$
Generating functions are a powerful tool for studying 
sequences and Wilf has written a wonderful text devoted entirely to their
study~\cite{wil:gen}.  There are several reasons why one might wish to
convert a sequence into its generating function.  It may be possible
to find a closed form for $f(x)$ when one does not exist for
$(a_n)_{n\ge0}$, or the expression for the generating function may be
used to derive one for the sequence.  Also, sometimes it is easier to obtain
various properties of the $a_n$, such as a
recursion or congruence relation, from $f(x)$ rather than directly.

By way of illustration, consider the sequence whose terms are
$$
p(n)=\mbox{ number of integer partitions of the number $n\in\bbZ_{\ge0}$}
$$
where an {\it integer partition\/} is a way of writing $n$ as an unordered sum
of positive integers.  For example, $p(4)=5$ because of the partitions
$$
4,\ 3+1,\ 2+2,\ 2+1+1,\ 1+1+1+1+1.
$$
There is no known closed form for $p(n)$, but the generating function
was found by Euler
\beq							\label{pgf}
\sum_{n=0}^\infty p(n)x^n=\prod_{k=1}^\infty \frac{1}{1-x^k}
\eeq
since $1/(1-x^k)=1+x^k+x^{2k}+\cdots$ and so a term in the product 
obtained by choosing $x^{ik}$ from this expansion corresponds to 
choosing a partition with the integer $k$ repeated $i$ times.
From~(\ref{pgf}) one can obtain all sorts of information about $p(n)$,
such as its asymptotic behavior.  See Andrews' book~\cite{and:tp} for
more details.

Our main object of study will be the generating function for the
M\"obius function of a poset $P$, the so-called characteristic 
polynomial.
Let $P$ have a zero and be {\it ranked\/} so that for any $x\in P$ all maximal 
chains from $\zh$ to $x$ have the same length denoted $\rho(x)$ and
called the {\it rank of $x$}.  (A
{\it chain\/} is a totally ordered subset of $P$ and {\it maximal\/}
refers to inclusion.)   The {\it
characteristic polynomial of $P$\/} is then
	\beq							\label{chi1}
	\chi(P,t)=\sum_{x\in P} \mu(x) t^{\rho(\oh)-\rho(x)}.
	\eeq
Note that we use the corank rather than the rank in the exponent on
$t$ so that $\chi$ will be monic.

\bfi
\bpi(70,90)(-10,-20)
\Gdb
\put(30,0){\makebox(0,0){$\{1\}\{2\}\{3\}$}}
\GadL{-50pt}{$\{1,2\}\{3\}$}
\GddL{-50pt}{$\{1,3\}\{2\}$}
\GgdL{-50pt}{$\{1\}\{2,3\}$}
\Gdf
\put(30,60){\makebox(0,0){$\{1,2,3\}$}}
\Gdbad \Gdbdd \Gdbgd \Gaddf \Gdddf \Ggddf
\epi
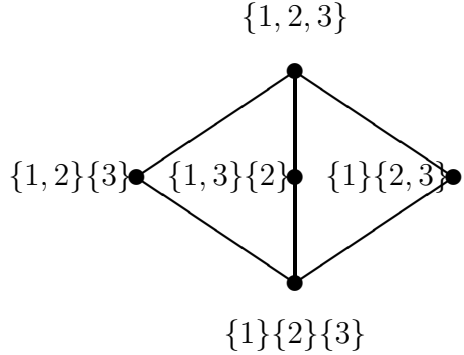
\capt{The partition poset $\Pi_3$}				\label{ex2}
\efi

Let us look at some examples of posets and their
characteristic polynomials, starting with those from the previous section.  
For the chain we clearly have
	$$\chi(C_n,t)=t^n-t^{n-1}=t^{n-1}(t-1).$$
Now for the Boolean algebra we have
	$$\chi(B_n,t)=\sum_{x\sbe[n]}(-1)^{|x|}t^{n-|x|}=(t-1)^n.$$
Note that by the same argument, if $k$ is the number of distinct
primes dividing $n$  then
$$
\chi(D_n,t)=t^{\rho(n)-k} (t-1)^k
$$
since the terms divisible by squares contribute nothing to the sum.
As a fourth example, consider the {\it partition poset $\Pi_n$\/}
which consists of all {\it set} partitions of $[n]$ (families of
disjoint nonempty subsets whose 
union is $[n]$)  ordered by refinement.  Direct computation with
$\Pi_3$ as shown in Figure~\ref{ex2} shows that
$\chi(\Pi_3,t)=t^2-3t+2=(t-1)(t-2)$.  In general
	$$\chi(\Pi_n,t)=(t-1)(t-2)\cdots(t-n+1).$$
Note that in all cases $\chi$ has only nonnegative integral roots.

Many of our example posets will arise as intersection lattices of
subspace arrangements.   A {\it lattice}, $L$, is poset such that 
every pair  $x,y\in L$ has a meet or greatest lower bound, $x\mt y$,
and a join or least upper bound, $x\jn y$.  All our lattices will be
finite and so will automatically have a zero $\zh=\Mt L$ and a one
$\oh=\Jn L$.  A {\it subspace arrangement\/} is a finite set
\beq							\label{cA}
	\cA=\{K_1,K_2,\ldots,K_l\}
\eeq
of subspaces of real Euclidean space $\bbR^n$.  If $\dim K_i=n-1$ for
$1\le i\le l$ then we say that $\cA$ is a {\it hyperplane arrangement}
and use $H$'s in place of $K$'s.
The {\it intersection lattice of $\cA$},  $L(\cA)$, has as elements
all subspaces $X\sbe\bbR^n$ that can be written as an intersection of
some of the elements of $\cA$.  The partial order is reverse
inclusion, so that $X\le Y$ if and only if $X\spe Y$.  So $L(\cA)$ has
minimal element $\bbR^n$, maximal element $K_1\cap\cdots\cap K_l$, and
join operation $X\jn Y=X\cap Y$.  The reader can
consult~\cite{bjo:sa,ot:ah} for more details about the general theory
of arrangements which is currently a very active field.

The {\it characteristic polynomial
of $A$\/} is defined by
	\beq							\label{chi2}
	\chi(\cA,t)=\sum_{X\in L(\cA)} \mu(X) t^{\dim X}.
	\eeq
This is not necessarily the same as $\chi(L(\cA),t)$ as defined
in~(\ref{chi1}).  If $\cA$ is a hyperplane arrangement then the two
will be equal up to a factor of a power of $t$, so from the point of
view of having integral roots there is no difference.  In the general
subspace case~\ree{chi1} and~\ree{chi2} may be quite dissimilar
and often the latter turns out to factor at least partially over
$\bbZ_{\geq0}$ while the former 
does not.  In the arrangement case the roots of~\ree{chi2} are called
the {\it exponents of $\cA$} and denoted $\exp \cA$.  In fact when $\cA$ is 
the set of reflecting hyperplanes for a Weyl group $W$ then these roots
are just the usual exponents of $W$~\cite{ter:gef} which are
always nonnegative integers.   The reason that Weyl groups,
as opposed to more general Coxeter groups, have well-behaved
characteristic polynomials is that such groups stabilize an
appropriate discrete subgroup of $\bbZ^n$.

All of our previous example lattices can be realized
as intersection lattices of subspace arrangements.  The $n$-chain is
$L(\cA)$ with $\cA=\{K_0,\ldots,K_n\}$ where $K_i$ is the set of all
points having the first $i$ coordinates zero.  The Boolean algebra is
the intersection lattice of the arrangement of coordinate hyperplanes
$H_i:\ x_i=0$, $1\le i\le n$.  By combining these two constructions,
one can also obtain a realization of the divisor poset as a subspace
arrangement.   To get the partition lattice we use the
{\it Weyl arrangement of type $A$}
$$
\cA_n=\{x_i-x_j=0\ :\ 1\le i<j\le n\}.
$$
To see why $\Pi_n$ and $L(\cA_n)$ are the same, associate the
hyperplane $x_i=x_j$ with the partition where $i,j$ are in one
subset and all other subsets are singletons.  This will then make the
join operations in the two lattices correspond.
Note that the characteristic polynomials defined by~(\ref{chi1})
and~(\ref{chi2}) 
are the same in the first two examples while $\chi(\cA_n,t)=t\chi(\Pi_n,t)$.

We will also be concerned with the hyperplane arrangements
associated with other Weyl groups.  The reader interested in more
information about these groups should consult the excellent text of
Humphreys~\cite{hum:rgc}.  In particular, the other two infinite
families 
\bea
\cB_n&=&\{x_i\pm x_j=0\ :\ 1\le i<j\le n\}\cup\{x_i=0\ :\ 1\le i\le n\},\\
\cD_n&=&\{x_i\pm x_j=0\ :\ 1\le i<j\le n\}
\eea
will play a role.  The corresponding characteristic polynomials are
listed in Table~\ref{weyl} along with $\chi(\cA_n,t)$ for
completeness.
(We do not consider the arrangement for the root system $C_n$ because
its roots are scalar multiples of the ones for $B_n$, thus yielding
the same arrangement.)
We will show how to derive the formulas for the characteristic
polynomials of $\cA_n,\cB_n$
and $\cD_n$ using elementary graph theory in the next section.

\begin{table}
$$\barr{l|l|l}
\cA	&\chi(\cA,t)			&\exp \cA\\ 
\hline
\cA_n	&t(t-1)(t-2)\cdots(t-n+1)	&0,1,2,\ldots,n-1\\
\cB_n   &(t-1)(t-3)\cdots(t-2n+1)	&1,3,5,\ldots,2n-1\\		
\cD_n	&(t-1)(t-3)\cdots(t-2n+3)(t-n+1)&1,3,5,\ldots,2n-3,n-1
\earr$$
\caption{Characteristic polynomials and exponents of some Weyl arrangements}
\label{weyl}
\end{table}

\section{Signed graphs}						\label{sg}

Zaslavsky developed his theory of signed
graphs~\cite{zas:grs,zas:sgc,zas:cis} to study hyperplane arrangements
contained in the Weyl arrangement $\cB_n$.  (Note that this includes
$\cA_n$ and $\cD_n$.)  In particular his coloring arguments provide
one of the simplest ways to compute the corresponding characteristic
polynomials.  

A {\it signed graph}, $G=(V,E)$, consists of a set $V$ of vertices
which we will always take to be $\{1,2,\ldots,n\}$, and a set of edges
$E$ which can be of three possible types:
\ben
\item a {\it positive edge\/} between $i,j\in V$, denoted $ij^+$,
\item a {\it negative edge\/} between $i,j\in V$, denoted $ij^-$,
\item and a {\it half-edge\/} which has only one endpoint $i\in V$,
denoted $i^h$. 
\een
One can have both the positive and negative edges between a given pair
of vertices in which case it is called a {\it doubled edge\/} and
denoted $ij^{\pm}$.
The three types of edges correspond to the three types of hyperplanes
in $\cB_n$, namely $x_i=x_j,\ x_i=-x_j,$ and $x_i=0$ for the positive,
negative, and half-edges, respectively.  
So to every hyperplane arrangement $\cA\sbe\cB_n$ we have an
associated signed graph $G_{\cA}$ with a hyperplane in $\cA$ if and
only if the corresponding edge is in $G_{\cA}$.
Actually, the possible edges in $G_{\cA}$ 
really correspond to the vectors in the root system of type $B_n$
perpendicular to the hyperplanes which are $e_i-e_j,\ e_i+e_j,$ and
$e_i$.  
(In the full theory one also considers the root system $C_n$
with roots $2e_i$ which are modeled by loops $ii$ in $G$. 
This is why the somewhat strange definition of a half-edge is
necessary.  Loops and half-edges behave differently because, e.g.,
the former can be in a circuit of the graph while the later can not.)
In picturing a signed graph I will  draw an ordinary edge for $ij^+$,
an edge with a slash through it for $ij^-$,
an edge with  two slashes through it for $ij^\pm$,
and an edge starting at a vertex and wandering off into space for $i^h$.
The graphs $G_{\cA_3},\ G_{\cB_3}$, and $G_{\cD_3}$ are
shown in Figure~\ref{graphs}.

\thicklines
\setlength{\unitlength}{1.7pt}
\bfi
\btab{ccc}
\bpi(70,70)(-10,-10)
\GaaL{-15pt}{2}
\GgaL{15pt}{3}
\GdgL{-15pt}{1}
\Gaaga \Gaadg \Ggadg
\put(30,-10){\makebox(0,0){$G_{\cA_3}$}}
\epi
&
\bpi(70,70)(-10,-10)
\GaaL{-15pt}{2}
\GgaL{15pt}{3}
\GdgL{-15pt}{1}
\Gaaga \Gaadg \Ggadg
\put(30,-10){\makebox(0,0){$G_{\cB_3}$}}
\put(27,-4){\line(0,1){8}}
\put(33,-4){\line(0,1){8}}
\put(10,30){\line(2,-1){8}}
\put(20,30){\line(-2,1){8}}
\put(40,30){\line(2,1){8}}
\put(50,30){\line(-2,-1){8}}
\put(0,0){\line(-1,-1){10}}
\put(60,0){\line(1,-1){10}}
\put(30,60){\line(0,1){10}}
\epi
&
\bpi(70,70)(-10,-10)
\GaaL{-15pt}{2}
\GgaL{15pt}{3}
\GdgL{-15pt}{1}
\Gaaga \Gaadg \Ggadg
\put(30,-10){\makebox(0,0){$G_{\cD_3}$}}
\put(27,-4){\line(0,1){8}}
\put(33,-4){\line(0,1){8}}
\put(10,30){\line(2,-1){8}}
\put(20,30){\line(-2,1){8}}
\put(40,30){\line(2,1){8}}
\put(50,30){\line(-2,-1){8}}
\epi
\etab
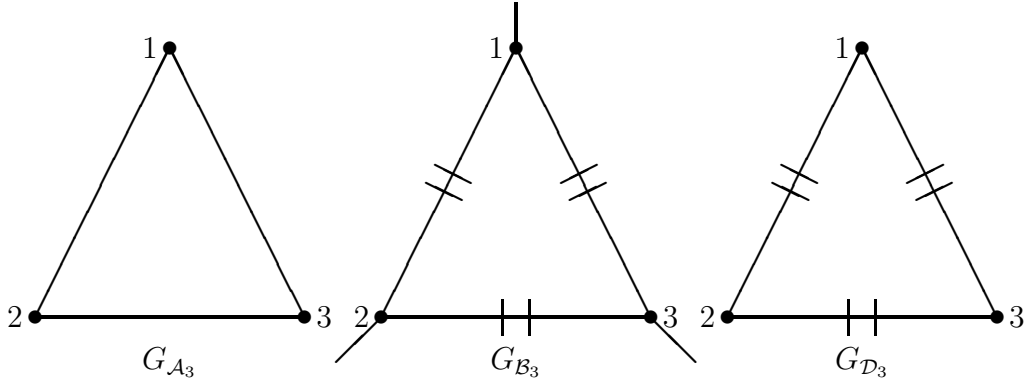
\capt{Graphs for Weyl arrangements}  \label{graphs}
\efi

Since we are using signed edges, we are also going
to use signed 
colors for the vertices.  For $s\in\bbZ_{\ge0}$ let
$[-s,s]=\{-s,-s+1,\ldots,s-1,s\}$.  A {\it coloring\/} of the signed
graph $G$ is a function $c:V\ra[-s,s]$.  The fact that the number of
colors $t=|[-s,s]|=2s+1$ is always odd will be of significance later.
A {\it proper coloring\/} $c$ of $G$ 
requires that  for every edge $e\in E$ we have
\ben
\item if $e=ij^+$ then $c(i)\neq c(j)$,
\item if $e=ij^-$ then $c(i)\neq -c(j)$,
\item if $e=i^h$ then $c(i)\neq 0$.
\een
Note that the first of these three restrictions is the one associated with
ordinary graphs and the four-color theorem~\cite{cl:gd}.  The {\it
chromatic polynomial of $G$\/} is 
$$
P(G,t)=\mbox{the number of proper colorings of $G$ with $t$ colors.}
$$
It is not obvious from the definition that $P(G,t)$ is actually a
polynomial in $t$.  In fact even more is true as we see in the
following theorem of Zaslavsky.  
\bth[\cite{zas:sgc}]  					\label{zas}
Suppose $\cA\sbe\cB_n$ has signed graph
$G_{\cA}$.  Then 
	$$\chi(\cA,t)=P(G_{\cA},t).\qquad\qed$$
\eth

Theorem~\ref{zas} trivializes the calculation of the characteristic
polynomials for the three infinite families of Weyl arrangements and
in so doing explains why they factor over $\bbZ_{\ge0}$.  For
$\cA_n$ the graph $G_{\cA_n}$ consists of every possible positive
edge.  So to properly color $G_{\cA_n}$ we have $t$ choices for
vertex 1, then $t-1$ for vertex 2 since $c(2)\neq c(1)$, and so forth
yielding
$$
\chi(\cA_n,t)=P(G_{\cA_n})=t(t-1)\cdots(t-n+1)
$$
in agreement with Table~\ref{weyl}.
It will be convenient in a bit to have a shorthand for this falling factorial,
so let $\fall{t}{n}=t(t-1)\cdots(t-n+1)$.
In $G_{\cB_n}$ we also have every negative edge and half-edge.  This
gives $t-1$ choices for vertex 1 since color 0 is not allowed, $t-3$
choices for vertex 2 since $c(2)\neq\pm c(1),0$, and so on.  These are exactly
the factors in the $\cB_n$ entry of Table~\ref{weyl}.  Finally
consider $G_{\cD_n}$ which is just 
$G_{\cB_n}$ with the half-edges removed.    There are two
cases depending on whether the color 0 appears once or not at all.
(It can't appear two or more times because $G_{\cA_n}\sbe G_{\cD_n}$.)
If the color 0 is never used then we have the same number of
colorings as with $\cB_n$.  If 0 is used once then there are $n$
vertices that could receive it and the rest are colored as in
$\cB_{n-1}$.  So
$$
\chi(\cD_n,t)=\prod_{i=1}^n(t-2i+1)+n\prod_{i=1}^{n-1}(t-2i+1)=
	(t-n+1)\prod_{i=1}^{n-1}(t-2i+1)
$$
which again agrees with the table.

Blass and I have generalized Zaslavsky's theorem from
hyperplane arrangements to subspace arrangements.  If $\cA$ and $\cB$
are subspace arrangements then we say that {\it $\cA$ is embedded in
$\cB$\/} if all subspaces of $\cA$ are intersections of subspaces of
$\cB$, i.e., $\cA\sbe L(\cB)$.  Now consider $[-s,s]^n$ as a cube of
integer lattice points in $\cR^n$  (not to be confused with our use of
lattice as a type of partially ordered set).  Let $[-s,s]^n\setm\bigcup\cA$
denote the set of points of the cube which lie on none of the
subspaces in $\cA$.  We will include a proof of the next result
because it amply illustrates the importance of the M\"obius Inversion
Theorem.
\bth[\cite{bs:cep}]						\label{cep}
Let $t=2s+1$ where $s\in\bbZ_{\ge0}$ and let $\cA$ be any subspace
arrangement embedded in $\cB_n$.  Then
	$$\chi(\cA,t)=|[-s,s]^n\setm\bigcup\cA|.$$
\eth
\pf\
 We construct two functions $f,g:L(\cA)\longrightarrow\bZ$ by
defining for each $X\in L(\cA)$
        \bea
        f(X)&=&|X\cap[-s,s]^n|,\\
        g(X)&=&|(X\setm\bigcup_{Y>X}Y)\cap[-s,s]^n|.
        \eea
Recall that $L(\cA)$ is ordered by {\it reverse} inclusion so that
$\bigcup_{Y>X}Y\sbs X$.  In particular $g(\bbR^n)=|[-s,s]^n\setm\bigcup\cA|$.   
Note also that $X\cap[-s,s]^n$ is combinatorially
just a cube of dimension $\dim X$ and side $t$ so that $f(X)=t^{\dim X}$.
Finally, directly from our definitions, $f(X) =\sum_{Y\ge X}g(Y)$ so by the 
Theorem~\ref{mit}
        \bea
        |[-s,s]^n\setm\bigcup\cA|&=&g(\zh)\\
                &=&\sum_{X\ge\zh)}\mu(X)f(X)\\
                &=&\sum_{X\in L(\cA)}\mu(X)t^{\dim X}\\
                &=&\chi(\cA,t)
        \eea
which is the desired result.\hfill\Qed

To see why our theorem implies Zaslavsky's in the hyperplane case,
note that a point $c\in[-s,s]^n$ is just a coloring
$c:V\ra[-s,s]$ where the $i$th coordinate of the point is the color of
the vertex $i$.  With this viewpoint, a coloring is proper if and only
if the corresponding point is not on any hyperplane of $\cA$.  For
example, if $ij^+\in E$ then the coloring must have $c(i)\neq c(j)$
and so the point does not lie on the hyperplane $x_i=x_j$.

As an application of Theorem~\ref{cep}, we consider a set of subspace
arrangements that has been arousing a lot of interest lately.  The
{\it $k$-equal arrangement of type $A$\/} is
$$
\cA_{n,k}=\{x_{i_1}=x_{i_2}=\ldots=x_{i_k}\ :\ 1\le i_1<i_2<\ldots<i_k\le n\}.
$$
The $\cA_{n,k}$ arrangements were introduced in the work
of Bj\"orner, Lov\'asz and Yao~\cite{bly:ldt} motivated,
surprisingly enough, by its
relevance to a certain problem in computational complexity.
Its study has been continued by many people including Linusson,
Sundaram, Wachs and
Welker~\cite{bl:ldt,bw:hkm,bw:nsc1,bw:nsc2,lin:mfc,lin:prb,sw:hrk,sw:gal}. 
Analogs of this subspace arrangement for types $B$ and $D$ have also
been studied by Bj\"orner and myself~\cite{bs:sat}.

Now in general $\chi(\cA_{n,k})$ does not factor completely over
$\bbZ_{\ge0}$, but it does factor partially.  In fact it is divisible
by the characteristic polynomial $\chi(\cA_m,t)=\fall{t}{m}$ for a certain
$m$.  What's more if one expands
$\chi(\cA_{n,k})$ in the basis $\fall{t}{n},$ $n\ge0$, for the
polynomial ring then the coefficients are nonnegative integers with a
simple combinatorial interpretation.
In particular,
let $S_k(n,j)$ denote the number of partitions of an $n$-element
set into $j$ subsets each of which is of size at most $k$.  Thus these
are generalizations of the Stirling numbers of the second kind.
We now have the expansion, first derived by Sundaram~\cite{sun:aht}
	\beq				\label{Seq}
	\chi(\cA_{n,k},t)=\sum_j S_{k-1}(n,j) \fall{t}{j}.
	\eeq
To see why this is true, consider an arbitrary point
$c\in[-s,s]^n\setm\bigcup\cA_{n,k}$.  So $c$ can have at
most $k-1$ of its coordinates equal.  Consider the $c$'s with exactly
$j$ different coordinates.  Then there are $S_{k-1}(n,j)$ ways to
distribute the $j$ values among the $n$ coordinates with at most $k-1$
equal.  Next we can choose which values to use in
$\fall{t}{j}$ ways.  Summing over all $j$ gives the desired equation.
From~(\ref{Seq}) we can immediately derive a divisibility relation.
To state it, let $\ce{\cdot}$ be the ceiling or round up function.  Then
$$
\fall{t}{\ce{n/(k-1)}}\ |\ \chi(\cA_{n,k},t)
$$
since $S_{k-1}(n,j)=0$ if
$j<\ce{n/(k-1)}$.   (Obviously $j$ sets of at most $k-1$ elements can
partition a set of size of at most $n=j(k-1)$.)

Thinking about things in terms of lattice points also permits a
generalization of Zaslavsky's theorem in another direction, namely to all
Weyl hyperplane arrangements (even the exceptional ones).  Let
$\Phi\sbs\bbR^n$ be a
root system for a finite Weyl group $W$ and let $\cW$ be the set of
hyperplanes perpendicular to the roots.  Let $(\cdot,\cdot)$ denote
the standard inner product on $\bbR^n$.  The role
of the cube in Theorem~\ref{cep} will be played by 
$$
P_t(\Phi)=\{x\in\bbR^n\ :\ \mbox{$(x,\al)\in\bbZ_{<t}$ for all $\al\in\Phi$}\}
$$
which is a set of points in the coweight  lattice  
of $\Phi$ closely associated
with the Weyl chambers of the corresponding affine Weyl group.

Consider a fixed simple system 
$$
\Delta=\{\sigma_1,\dots,\sigma_n\}
$$ 
in $\Phi$.  Since $\Delta$ is a basis for
$\bR^n$ any root $\al\in\Phi$ can be
written as a linear combination,
$$
\al=\sum_{i=1}^n s_i(\al)\sigma_i.
$$
In fact the coefficients $s_i(\al)$ are always integers.
Among all the roots, there is a {\it highest} one, $\alt$,
characterized by the fact that for all roots $\al$ and all
$i\in[n]$, $s_i(\alt)\geq s_i(\al)$.  
We will also need a weighting factor called the {\it
index of connection}, $f$, which is the index of 
the lattice generated by the roots in the coweight lattice.
The second generalization can now be stated.
\bth[\cite{ath:cps,bs:cep,hai:cqr}]				\label{cep2}
Let $\Phi$ be a root system for a finite Weyl group with associated
arrangement $\cW$.  Let $t$ be a positive integer relatively prime
to $s_i(\alt)$ for all $i$.  Then 
\beq							\label{cepeq}
\chi(\cW,t)=\frac{1}{f}\left|P_t(\Phi)\setm\bigcup\cW\right|.\qquad\qed
\eeq
\eth
Note how the condition in Theorem~\ref{cep} that $t$ be odd has
been replaced by a relative primeness restriction.  This is typical when
dealing with Ehrhart quasi-polynomials~\cite[page 235ff.]{sta:ec1} 
which enumerate the number of
points of a given lattice inside a polytope and its blowups.
We have not been able to
use~(\ref{cepeq}) to explain the factorization of $\chi(\cW,t)$ over
$\bZ_{\ge0}$ as was done with Theorem~\ref{zas} for the three infinite
families.  It would be interesting if this hole could be filled.

Athanasiadis~\cite{ath:cps} has given a very pretty proof of the
previous theorem.
His main tool is a reworking of a result of Crapo and Rota~\cite{cr:cg}
which is similar 
in statement and proof to Theorem~\ref{cep} but
replaces $[-s,s]^n$ by $\bbF_p^n$ where $\bbF_p$ is the finite
field with $p$ elements, $p$ prime.  Terao~\cite{ter:jdf} also independently
discovered this theorem.
\bth[\cite{ath:cps,cr:cg,ter:jdf}]			\label{cr}
Let $\cA$ be any subspace arrangement in
$\bbR^n$ defined over the integers and hence over $\bbF_p$.  Then for
large enough primes $p$ we have
$$
\chi(\cA,p)=|\bbF_p^n\setm\bigcup\cA|.\quad\qed
$$
\eth

Athanasiadis has also used the previous result to give very elegant
derivations of characteristic polynomials for many arrangements which
cannot be handled by Theorem~\ref{cep}.  For the rest of this section
only, we will enlarge the definition of a hyperplane arrangement to be
any finite set of affine hyperplanes (not necessarily passing through
the origin).  An arrangement $\cB$ is a {\it deformation} of arrangement 
$\cA$ if every hyperplane of $\cB$ is parallel to some hyperplane of 
$\cA$.  As an example consider the {\it Shi arrangement of type A}, $\cS_n$,
with hyperplanes
$$
x_i-x_j=0,\ x_i-x_j=1,\qmq{where} 1\le i<j\le n	
$$
which is a deformation of the corresponding Weyl arrangement.
Such arrangements were introduced by Shi~\cite{shi:klc,shi:stc}
for studying affine Weyl groups.  Headly~\cite{hea:rei,hea:rwa} first
computed their characteristic polynomials in a way that
was universal for all types but relied on a formula of Shi's that had
a complicated proof.  To illustrate the power of Theorem~\ref{cr}, we
will reproduce the proof in~\cite{ath:cps} that
\beq						\label{chiS}
\chi(\cS_n,t)=t(t-n)^{n-1}.
\eeq
Consider any $(x_1,\ldots,x_n)\in\bbF_p^n$
as a placement of balls labeled $1,\ldots,n$ into a circular array of
boxes labeled clockwise as $0,\ldots,p-1$, where placement of ball $i$
in box $j$ means $x_i=j$.  Then $(x_1,\ldots,x_n)\in\bbF_p^n\setm\bigcup\cS_n$
means that no two balls are in the same box and that if two balls are
in adjacent boxes then they must be in increasing order clockwise.
All such placements can be derived as follows.  Take $p-n$ unlabeled
boxes and place them in a circle.  Now put the balls $1,\ldots,n$
into the spaces between the boxes so that adjacent ones
increase clockwise.  Note
that since the boxes 
are unlabeled, there is only one way to place ball $1$, but once that is
done balls $2,\ldots,n$ can be placed in $(p-n)^{n-1}$ ways.  Finally, put
an unlabeled box around each ball, label all
the boxes clockwise as $0,\ldots,p-1$ which can be done in $p$ ways,
and we are done.

The connected components of $\bbR^n\setm\cA$ are called {\it regions}
and a region is {\it bounded} if it is contained in some sphere about
the origin.  If we let $r(\cA)$ and $b(\cA)$ denote the number of regions
and bounded regions, respectively, of $\cA$ then we can state the
following striking result of Zaslavsky.
\bth[\cite{zas:fua}]
For any affine hyperplane arrangement 
$$
r(\cA)=|\chi(\cA,-1)|=\sum_{X\in L(\cA)}|\mu(X)|
$$
and
$$
b(\cA)=|\chi(\cA,1)|=\left|\sum_{X\in L(\cA)}\mu(X)\right|.\quad\qed
$$
\eth
Using the characteristic polynomials in Table~\ref{weyl} we see that
$$
\barr{lclcl}
b(\cA_n)&=&|(-1)(-1-1)(-1-2)\cdots(-1-n+1)|	&=&n!\\
b(\cB_n)&=&|(-1-1)(-1-3)\cdots(-1-2n+1)|	&=&2^n n!\\
b(\cC_n)&=&|(-1-1)(-1-3)\cdots(-1-2n+3)(-1+n+1)|&=&2^{n-1}n!\\
\earr
$$
which agrees with the well-known fact that the number of chambers of a
Weyl arrangement is the same as the number of elements in the
corresponding group.
It was Shi's formula for the number of regions in his
arrangements that Headly needed to derive the full characteristic
polynomial.  In particular, it follows from~\ree{chiS} that
$r(\cS_n)=(n+1)^{n-1}$ 
which is known to be the number of labeled trees on $n+1$
vertices, or the number of parking functions of length $n$.  Many
people~\cite{al:sbr,hea:rwa,pos:eag,ps:dch,sta:hai} have have used these
combinatorial interpretations to come up with bijective proofs of this formula
and related ones.

\section{Free arrangements}				\label{fa}

In this section we consider a large class of hyperplane arrangements
called free arrangements which were introduced by
Terao~\cite{ter:gef}.  The characteristic polynomial of such an
arrangement factors over $\bbZ_{\ge0}$ because its roots are related to
the degrees of basis elements for a certain associated free module.

Our modules will be over the polynomial algebra
$A=\bbR[x_1,\ldots,x_n]=\bbR[x]$ 
graded by total degree so $A=\oplus_{i\geq0} A_i$.
The {\it module of derivations}, $D$, consists of all $\bbR$-linear maps
$\theta:A\to A$ satisfying
	$$\theta(fg)=f\theta(g)+g\theta(f)$$
for any $f,g\in A$.  
This module can be graded by saying that $\theta$ has degree $d$,
$\deg\th=d$,  if 
$\theta(A_i)\sbe A_{i+d}$ for all $i\ge0$.  Also, $D$ is free with  basis
$\partial/\partial x_1, \ldots, \partial/\partial x_n$.
It is simplest to display a derivation as a column vector
with entries being its components with respect to this basis.  So if
$\theta=p_1(x)\partial/\partial x_1 + \cdots +p_n(x)\partial/\partial x_n$
then we write
	$$\theta = \left[\barr{c} p_1(x) \\
				\vdots \\
				p_n(x)
		\earr
		\right]
		=  \left[\barr{c} \theta(x_1) \\
				\vdots \\
				\theta(x_n)
		\earr
		\right].$$
Two operators that we will find useful are
$$
\bX^d=x_1^d\partial/\partial x_1 + \cdots +x_n^d\partial/\partial x_n=
\left[\barr{c} x_1^d\\
		\vdots \\
		x_n^d \earr\right]
$$
and
$$
\bXh=\xh_1\partial/\partial x_1 + \cdots +\xh_n\partial/\partial x_n=
\left[\barr{c} \xh_1\\
		\vdots \\
		\xh_n \earr\right]
$$
where $\hat{x}_i=x_1 x_2 \cdots  x_n/x_i$.
Note that $\deg\bX^d=d-1$ and $\deg\bXh=n-2$.

To see the connection with hyperplane arrangements, notice that any
hyperplane $H$ is defined by a linear equation
$\al_{H}(x_1,\ldots,x_n)=0$.  It is then useful to study the  
associated {\it module of $\cA$-derivations} which is defined by
$$
D(\cA)=\{\th\in D\ :\ \mbox{$\al_H|\th(\al_H)$ for all $H\in\cA$}\}
$$
where $p|q$ is division of polynomials in $A$.  
By way of illustration, $\bX^d\in D(\cA_n)$ for all $d\ge0$
since $\bX^d(x_i-x_j)=x_i^d-x_j^d$ which is divisible by $x_i-x_j$.
Similarly $\bX^{2d+1}\in D(\cD_n)$ because of what we just showed for
$\cA_n$ and the fact that $x_i+x_j|x_i^{2d+1}+x_j^{2d+1}$.  The
$\bX^{2d+1}$ are also in $D(\cB_n)$ since $x_i|x_i^{2d+1}$.
By the same methods we get $\bXh\in D(\cD_n)$.

We say that the
arrangement $\cA$ is {\it free\/} if $D(\cA)$ is free as a module over $A$.
Freeness is intimately connected with the factorization of $\chi$ as
the next theorem shows.
\bth[\cite{ter:gef}]					\label{DA}
If $\cA$ is free then
$D(\cA)$ has a homogeneous basis $\theta_1,\ldots,\theta_n$ and
the degree set
$\{d_1,\ldots, d_n\}=\{\deg \theta_1, \ldots, \deg \theta_n\}$
depends only on $\cA$.  Furthermore
	$$\chi(\cA,t)=(t-d_1-1)\cdots(t-d_n-1).\qquad\qed$$
\eth

We have a simple way to check whether a derivation is in $D(\cA)$ for a
given arrangement $\cA$.  It would be nice to have an easy way to test
whether $\cA$ is free and if so find a basis.  This is the
Saito Criterion.
Given derivations $\theta_1,\ldots,\theta_n$,
consider the 
matrix whose columns are the corresponding column vectors
	$$\Th = [\theta_1, \ldots, \theta_n] = [\theta_j(x_i)].$$
Also consider the homogeneous polynomial
	$$Q=Q(\cA)=\prod_{H\in\cA} \al_H(x)$$
which has the arrangement $\cA$ as zero set.
For example
\bea
Q(\cA_n)&=&\prod_{1\le i<j\le n} (x_i-x_j)\\
Q(\cB_n)&=&x_1x_2\cdots x_n\prod_{1\le i<j\le n} (x_i^2-x_j^2)\\
Q(\cD_n)&=&\prod_{1\le i<j\le n} (x_i^2-x_j^2).
\eea
\bth[\cite{sai:tld,ter:fah}]				\label{det}
Suppose $\theta_1,\ldots,\theta_n\in D(\cA)$ and that $Q$ is the
defining form of $\cA$.   Then $\cA$ is free with basis
$\theta_1,\ldots,\theta_n$  if and only if
$$
\det\Th=cQ
$$ 
for some $c\in\bbR\setm0$.\hfill\qed
\eth

How could this be applied to the  Weyl arrangements?  Given what we know about
elements in their derivation modules and the factorization of their
characteristic polynomials, it is natural to guess that we might be
able to prove freeness with the following matrices
\bea
\Th(\cA_n)&=&\left[\bX^0, \bX^1, \bX^2, \ldots, \bX^{n-1}\right],\\
\Th(\cB_n)&=&\left[\bX^1, \bX^3, \bX^5, \ldots, \bX^{2n-1}\right],\\
\Th(\cD_n)&=&\left[\bX^1, \bX^3, \bX^5, \ldots, \bX^{2n-3},\bXh\right].
\eea
Of course $\det\Th(\cA_n)=\prod_{1\le i<j\le n} (x_i-x_j)=\pm Q(\cA_n)$
is just Vandermonde's determinant.   Similarly we get
$\det\Th(\cB_n)=\pm x_1x_2\cdots x_n\prod_{1\le i<j\le n}(x_i^2-x_j^2)$
by first factoring out $x_i$ from the $i$th row which results in a
Vandermonde in squared variables.  
For $\cD_n$ just factor out $x_1x_2\cdots x_n$ from the last column
and then put these factors back in by multiplying row $i$ by $x_i$.
The result is again a Vandermonde in squares.
Now the roots of the corresponding characteristic
polynomials can be read off these matrices in agreement with Table~\ref{weyl}.

The reader may have noticed that the bases we have for
$D(\cB_n)$ and $D(\cD_n)$ are
the same except for the last derivation.  This is reflected in the fact that
$\exp{\cB_n}$ and $\exp{\cD_n}$ are the same except for the last root.
Note that the difference between these roots is $n$ which is exactly
the number of hyperplanes in $\cB_n$ but not in $\cD_n$.  Wouldn't it
be lovely if adding these hyperplanes one at a time to $\cD_n$ would
produce a sequence of arrangements all of whose exponents agreed with
$\exp(\cD_n)$ except the last one which would increase by one each
time a hyperplane is added until we reach $\exp(\cB_n)$?  This is in
fact what happens.  Define
$$
\cDB_{n,k}=\cD_n\cup\{x_1,x_2,\ldots,x_k\}
$$
so that $\cDB_{n,0}=\cD_n$ and $\cDB_{n,n}=\cB_n$.  Now the derivation
$\th_k=x_1x_2\cdots x_k\bXh$ (scalar multiplication) is in $D(\cDB_{n,k})$
since $\bXh\in D(\cD_n)$ and $x_i\ |\ \th_k(x_i)$ for $1\le i\le k$.
Furthermore, if we let
$$
\Th(\cDB_{n,k})=\left[\bX^1, \bX^3, \bX^5, \ldots, \bX^{2n-3},\th_k\right]
$$
then $\det\Th(\cDB_{n,k})=x_1x_2\cdots x_k\det\Th(\cD_n)=Q(\cDB_{n,k})$ so 
we do indeed have a basis.  Thus 
$\exp(\cDB_{n,k})=\{1,3,5,\ldots,2n-3,n-1+k\}$ as desired.
The $\cDB_{n,k}$ were first considered by Zaslavsky~\cite{zas:grs}.
Bases for the module of derivations associated to other hyperplane
arrangements interpolating between the three infinite Weyl families
have been computed by J\'ozefiak and myself~\cite{js:bds}.  Edelman
and Reiner~\cite{er:fha} have determined all free arrangements lying
between $\cA_n$ and $\cB_n$.  It is still an open problem to find all
the free subarrangements of $\cB_n$ which do not contain $\cA_n$.

Related to these interpolations are the notions of inductive and
recursive freeness.  If $\cA$ is any hyperplane arrangement and
$H\in\cA$ then we have the corresponding {\it deleted arrangement\/}
$$
\cA'=\cA\setm H
$$
and the {\it restricted arrangement\/}
$$
\cA''=\{H'\cap H\ :\ H'\in\cA'\}.
$$
In this case $(\cA,\cA',\cA'')$ is called a {\it triple of arrangements}.
Of course $\cA'$ and $\cA''$ depend on $H$ even though the notation
does not reflect this fact.  Also if $\cA\sbe\cB_n$ then one can mirror
these two operations by defining deletion or contraction of
corresponding edges in
$G_{\cA}$.  The following Deletion-Restriction Theorem
shows how the characteristic polynomials for these three arrangements
are related.
\bth[\cite{bry:bcc,zas:fua}]
If $(\cA,\cA',\cA'')$ is a triple of arrangements then
$$
\chi(\cA,t)=\chi(\cA',t)-\chi(\cA'',t).\qquad\qed
$$
\eth
For freeness, we have Terao's Addition-Deletion
Theorem.  Note that its statement about the exponents
follows immediately from the previous result.
\bth[~\cite{ter:ahf}]					\label{ad}
If $(\cA,\cA',\cA'')$ is a triple of arrangements then any two of the
following statements implies the third:
\bea
\mbox{$\cA$ is free with} &&\exp\cA=\{e_1,\ldots,e_{n-1},e_n\},\\ 
\mbox{$\cA'$ is free with} &&\exp\cA'=\{e_1,\ldots,e_{n-1},e_n-1\},\\ 
\mbox{$\cA''$ is free with} &&\exp\cA''=\{e_1,\ldots,e_{n-1}\}.\qquad\qed
\eea
\eth

Continuing to follow~\cite{ter:ahf}, define the class $\cIF$
of {\it inductively free arrangements\/} to be those generated by the rules
\ben
\item[(1)] the empty arrangement in $\bbR^n$ is in $\cIF$ for all $n\ge0$,
\item[(2)] if there exists $H\in\cA$ such that $\cA',\cA''\in\cIF$ and
$\exp\cA''\sbs\exp\cA'$ then $\cA\in\cIF$.  
\een
So to show that $\cA$ is inductively free, we must start with an
arrangement which is known to be inductively free and add hyperplanes
one at a time so that~(2) is always satisfied.
If $\cF$ denotes the class of free arrangements then Theorem~\ref{ad}
shows that $\cIF\sbe\cF$ and one can come up with examples to show
that the inclusion is in fact strict.  On the other hand, it is not
hard to show using interpolating arrangements that $\cA_n, \cB_n$ and
$\cD_n$ are all inductively free.  Ziegler~\cite{zie:ach} has
introduced an even larger class of arrangements.  The class of {\it
recursively free arrangements}, $\cRF$, is gotten by using the same
two conditions as for $\cIF$ plus
\ben
\item[(3)] if there exists $H\in\cA$ such that $\cA,\cA''\in\cRF$ and
$\exp\cA''\sbs\exp\cA$ then $\cA'\in\cIF$.  
\een
It can be shown that $\cIF\sbs\cRF$ strictly but it is not known
whether every free arrangement is recursively free.

\section{Supersolvability}

In this section we will look at a combinatorial method
of Stanley~\cite{sta:sl} which  applies to lattices in general,
not just those which arise from arrangements.  First, however, we must
review an important result of Rota~\cite{rot:tmf} which gives a combinatorial
interpretation to the M\"obius function of a semimodular lattice.

A lattice $L$ is {\it modular\/} if for all $x,y,z\in L$ with $y\le z$
we have an associative law
\beq							\label{mod}
y\jn(x\mt z)=(y\jn x)\mt z.
\eeq
A number of natural examples, e.g., the partition lattice, are not
modular but satisfy the weaker condition
$$
\mbox{if $x$ and $y$ both cover $x\mt y$ then $x\jn y$ covers both $x$
and $y$}
$$
for all $x,y\in L$. (If $x,y\in L$ then {\it $x$ covers $y$\/} if
$x>y$ and there is no $z$ with $x>z>y$.)  Such lattices are called
{\it semimodular}. 
Lattice $L$ is modular if and only if both $L$ and its dual $L^*$
(where the order relation is reversed) are semimodular.

A set of important elements of $L$
are its {\it atoms\/} which are all elements $a$ covering $\zh$.  
We let $A(L)$ denote the atom set of $L$. 
If $L$ is semimodular then one can show that it is  ranked.
Furthermore, if $B\sbe A(L)$ then one can prove
that 
\beq							\label{ind}
\rho(\Jn B)\le |B|.
\eeq
We will call $B$ {\it independent\/} and a
{\it base\/} for $x=\Jn B$ if~(\ref{ind}) holds with equality.
This terminology comes from the theory of vector spaces.  Indeed if one
takes $L$ to be the lattice of all subspaces of $\bbF_q^n$ ($\bbF_q$ a
finite field) ordered by inclusion then atoms have dimension 1 and
lattice independence corresponds to independence of lines.  A
{\it circuit\/} is a dependent set which is minimal with respect to
inclusion.   If arrangement $\cA\sbe\cA_n$ has graph $G=G_{\cA}$ then
the atoms of $L(\cA)$ are edges of $G$ and a circuit of $L(\cA)$ 
forms a circuit in $G$ in the usual graph-theoretic sense.

\thicklines
\setlength{\unitlength}{1.3pt}
\bfi

\bpi(140,140)(-70,-70)
\put(0,-60){\circle*{3}}
\put(-60,-20){\circle*{3}}
\put(-20,-20){\circle*{3}}
\put(20,-20){\circle*{3}}
\put(60,-20){\circle*{3}}
\put(0,60){\circle*{3}}
\put(-60,20){\circle*{3}}
\put(-20,20){\circle*{3}}
\put(20,20){\circle*{3}}
\put(60,20){\circle*{3}}
\put(0,-70){$\zh$}
\put(-70,-20){$a$}
\put(-30,-20){$b$}
\put(30,-20){$c$}
\put(70,-20){$d$}
\put(0,70){$\oh$}
\put(-70,20){$s$}
\put(-30,20){$t$}
\put(30,20){$u$}
\put(70,20){$v$}
\put(0,-60){\line(-3,2){60}}
\put(0,-60){\line(-1,2){20}}
\put(0,-60){\line(1,2){20}}
\put(0,-60){\line(3,2){60}}
\put(0,60){\line(-3,-2){60}}
\put(0,60){\line(-1,-2){20}}
\put(0,60){\line(1,-2){20}}
\put(0,60){\line(3,-2){60}}
\put(-60,-20){\line(0,1){40}}
\put(-20,-20){\line(-1,1){40}}
\put(20,-20){\line(-2,1){80}}
\put(-60,-20){\line(1,1){40}}
\put(-20,-20){\line(1,1){40}}
\put(20,-20){\line(1,1){40}}
\put(60,-20){\line(-2,1){80}}
\put(60,-20){\line(-1,1){40}}
\put(60,-20){\line(0,1){40}}
\epi
\capt{A lattice $L$}				    		\label{L}
\efi

Now impose an arbitrary total order on $A(L)$ which will be denoted
$\lte$ so as to distinguish it from the partial order $\le$ on $L$.
A circuit $C\sbe A(L)$ gives rise to a {\it broken circuit\/}, $\Cb$,
obtained by removing the minimal element of $C$ in $\lte$.  A set
$B\sbe A(L)$ is {\it NBC\/} (No Broken Circuit) if $B$ does not
contain any of the $\Cb$.  Note that in this case $B$ must be
independent and so a base for $\Jn B$.  To illustrate, consider the
semimodular lattice $L$ in Figure~\ref{L}.  If we order the atoms
$a\lt b\lt c\lt d$ then $L$ has unique circuit $C=\{a,b,c\}$ with
associated broken circuit $\Cb=\{b,c\}$.  Comparing the number of NBC
bases of each element with its M\"obius function in the following
table
$$\barr{r|cccccccccc}
\mbox{element $x$}&\zh		&a	&b	&c	&d	&s
&t      &u	&v	&\oh \\ 
\hline
\mbox{NBC bases of $x$}&\emptyset&a	&b	&c	&d	&a,b
&a,d	&b,d	&c,d	&a,b,d \\
	        &	  	& 	& 	& 	& 	&a,c
&   	&   	&   	&a,c,d\\
\mu(x)			&+1	&-1	&-1	&-1	&-1	&+2
&+1	&+1	&+1	&-2
\earr$$
should lead the reader to a conjecture!  This is in fact the famous
result of Rota referred to earlier and usually called the NBC Theorem.
\bth[\cite{rot:tmf}]					\label{nbc}
Let $L$ be a semimodular lattice.  Then for any total ordering $\lte$
of $A(L)$ we have
$$
\mu(x)=(-1)^{\rho(x)}(\mbox{number of NBC bases of $x$}).\qquad\qed
$$
\eth

In order to apply the NBC theorem to our factorization problem, we
will need to make an additional restriction on $L$.  Write $xMz$  and
call $x,z$ a
{\it modular pair\/} if equation~(\ref{mod}) is satisfied for all $y\le z$.
Furthermore $x\in L$ is a {\it modular element\/} if $xMz$ and $zMx$
for every $z\in L$.  For example, if $L=L(\cA_n)$ or
$L(\cB_n)$ then an element corresponding to a graph $K_W$ which has a
complete component on the vertex set $W\sbe\{1,2,\ldots,n\}$ (all
possible edges from the parent graph between elements of $W$)  
and all other components trivial (isolated vertices) is modular.
A semimodular lattice is {\it supersolvable\/} if it has a maximal
chain of modular elements.  The lattice of subgroups of a finite
supersolvable group (one possessing a normal series where quotients of
consecutive terms are cyclic) is supersolvable.  From the previous
example we see that $L(\cA_n)$ and $L(\cB_n)$ are supersolvable.
However it is not true that $L(\cD_n)$ is supersolvable as we will see
later.  

Now any chain
$\zh=x_0<x_1<\ldots<x_n=\oh$ in $L$ defines a partition of the atoms
$A(L)$ into subsets
\beq								\label{Ak}
A_i=\{a\in A(L)\ :\ \mbox{$a\le x_i$ and $a\not\le x_{i-1}$}\}
\eeq
called {\it levels}.
A total order $\lte$ on $A(L)$ is said to be {\it induced\/} if it satisfies
\beq								\label{indu}
\mbox{if $a\in A_i$ and $b\in A_j$ with $i<j$ then $a\lte b$.}
\eeq
With these definitions we can state one of Stanley's main
results~\cite{sta:sl} about semimodular supersolvable lattices.
It states that their 
characteristic polynomials 
factor over $\bbZ_{\ge0}$ because the roots are the cardinalities of
the the $A_i$.
\bth[\cite{sta:sl}]  						\label{sl}
Let $L$ be a semimodular supersolvable lattice and suppose
$\zh=x_0<x_1<\ldots<x_n=\oh$ is a maximal chain of modular elements of $L$.
Then for any induced total order $\lte$ on $A(L)$
\ben
\item[(1)] the NBC bases of $L$ are exactly the sets of atoms gotten by
picking at most one atom from each $A_i$,
\item[(2)] $\chi(L,t)=(t-|A_1|)(t-|A_2|)\cdots(t-|A_n|).$
\een
\eth
\pf\  The proof that (1) implies (2) is so simple and beautiful that
I cannot resist giving it.  The coefficient of $t^{n-k}$ on the right
side of (2) is $(-1)^k$ times the number of ways to pick atoms from
exactly $k$ of 
the $A_i$.  But by (1) this is up to sign the number of NBC bases of
elements at rank $k$.  Putting back in the sign and using the NBC
theorem, we see that this coefficient is the sum of all the M\"obius
values for elements of rank $k$, which agrees with the corresponding
coefficient on the left side.\hfill\qed

As an example, consider the chain of graphs with a single nontrivial
complete component
$$
\zh=K_{\{1\}}<K_{\{1,2\}}<\ldots<K_{\{1,2,\ldots,n\}}=\oh
$$
in $\Pi_n=L(\cA_n)$.  Then $A_k$ is the set of all positive edges from $k+1$
to $i$, $i<k+1$, and so $|A_k|=k$.   Thus
$\chi(\Pi_n,t)=\prod_{i=1}^{n-1}(t-i)$ as before.  Using the
analogous chain in $L(\cB_n)$ (which starts at $K_\emp$) gives $A_k$
as containing all edges $ik^{\pm}, i<k$, and all half-edges $j^h, j\le k$.
So $|A_k|=2k-1$ giving the usual roots.  Now we can also see why
$L(\cD_n)$ is not supersolvable for $n\ge4$.  When $n\ge4$ the second
smallest root of $\chi(\cD_n,t)$ is 3.  So if the lattice were
supersolvable then Theorem~\ref{sl} would imply that some element
$x\in L(\cD_n)$ of rank two would have to cover at least $3+|A_1|=4$
atoms.  It is easy to verify that there is no such element.

It is frustrating that $L(\cD_n)$ is not supersolvable.  To get around
this problem, Bennett and I have introduced a more general
concept~\cite{bs:gss}.  Looking at the previous proof, the reader will
note that it would still go through if every NBC base could be obtained
in the following manner.  First pick an atom from a set
$A_1=\{a_1,a_1',a_1'',\ldots\}$.
Then pick the second atom from one of a family of sets
$A_2,A_2',A_2'',\ldots$ according to whether the first atom picked was
$a_1,a_1',a_1'',\ldots$ respectively, where $|A_2|=|A_2'|=|A_2''|=\ldots$,
and continue similarly.  This process can be modeled by an object
which we call an {\it atom decision tree\/} or {\it ADT\/} and any
lattice admitting an ADT has a characteristic polynomial with roots
$r_i$ equal to the common cardinality of all the sets of index $i$.
It turns out that the lattices for all of the interpolating
arrangements $\cDB_{n,k}$ admit ADTs and this combinatorially explains
their factorization.  H\'el\`ene Barcelo and Alain
Goupil~\cite{bg:nbc} have independently
come up with a factorization of the NBC complex of $L(\cD_n)$ (the
simplicial complex of all NBC bases of a lattice)
which is similar to the ADT one.   Their paper also contains a
nice result (joint with Garsia) relating the NBC
sets with reduced decompositions into reflections of Weyl group
elements.

Another way to generalize the previous theorem is to replace the
semimodularity and supersolvability restrictions by weaker conditions.
The new concepts are based on a generalization of the NBC Theorem that
completely eliminates semimodularity from its hypothesis.  Let $\lte$
be any {\it partial\/} order on $A(L)$.  It can be anything from a
total order to the total incomparability order induced by the ordering
on $L$.  A set $D\sbe A(L)$ is {\it bounded below\/} if for any 
$d\in D$ there is $a\in A(L)$ such that
\ben
\item[(a)] $a\lt d$ and
\item[(b)] $a<\Jn D$.
\een
In other words $a$ bounds $d$ below in $(A(L),\lte)$ and also bounds
$\Jn D$ below in $(L,\le)$.  We say $B\sbe A(L)$ is {\it NBB\/} if it
contains no bounded below set and say that $B$ is an {\it NBB base\/}
for $x=\Jn B$.  Blass and I have proved the following NBB Theorem
which holds for any lattice.
\bth[\cite{bs:mfl}]  Let $L$ be any lattice and let $\lte$ be any
partial order on $A(L)$.  Then for any $x\in L$ we have
$$
\mu(x)=\sum_{B} (-1)^{|B|}
$$
where the sum is over all NBB bases of $x$.\hfill\qed
\eth
Note that when $L$ is semimodular and $\lte$ is total then the NBB and
NBC bases coincide.  Also in this case all NBC bases of $x$ have the
same cardinality, namely $\rho(x)$, and so our theorem reduces to
Rota's.  However this result has much wider applicability, giving
combinatorial explanations for the M\"obius functions of the
non-crossing partition lattices and their type $B$ and $D$
analogs~~\cite{kre:pnc,rei:npc},
integer partitions under dominance
order~\cite{bog:mfd,bry:lip,gre:clm}, and 
the shuffle posets of Greene~\cite{gre:ps}.

Call $x\in L$ {\it left modular\/} if $xMz$ for all $z\in L$.  So this
is only half of the condition for modularity of $x$.  Call $L$ itself
{\it left modular\/} if
\bce
$L$ has a maximal chain
$\zh=x_0<x_1<\ldots<x_n=\oh$ of left modular elements.
\ece
This is strictly weaker than supersolvability as can be seen by
considering the 5-element nonmodular 
lattice~\cite[Proposition 2.2 and ff.]{sta:sl}.

In Stanley's theorem we cannot completely do away with semimodularity
as we did in Rota's (the reason why will come shortly),  but we can
replace it with a weaker hypothesis which we call the
{\it level condition}.  In it we assume that the partial order $\lte$
has been induced by some maximal chain, i.e., satisfies~(\ref{indu})
with ``if'' replaced by ``if and only if.''
\begin{center}
If $\lte$ is induced and $b_0\lt b_1\lt b_2\lt\dots\lt b_k$ then
$b_0\not\leq\Jn_{i=1}^kb_i$.  
\end{center}
It can be shown that semimodularity implies the level condition for
any induced order but
not conversely.  An {\it LL lattice\/} is one having a maximal left
modular chain such that the induced partial order satisfies the level
condition. So Theorem~\ref{sl} generalizes to the following.  Note that we must
extend the definition of the characteristic polynomial since an LL
lattice may not have a rank function and the first of the two parts makes
$\chi$ well-defined.
\bth[\cite{bs:mfl}]
Let $L$ be an an LL lattice with $\lte$ the 
partial order on $A(L)$ induced by a left modular chain.
\ben
\item[(1)] 
The NBB bases of $L$ are exactly the sets of atoms obtained by
picking at most one atom from each $A_i$ and all NBB  bases of a
given $x\in L$ have the same cardinality denoted $\rho(x)$.
\item[(2)] If we define $\chi(L,t)=\sum_{x\in L} \mu(x) t^{\rho(\oh)-\rho(x)}$
with $\rho$ as in (1), then
$$
\chi(L,t)=(t-|A_1|)(t-|A_2|)\cdots(t-|A_n|).\qquad\qed
$$
\een
\eth
This theorem can be used on lattices where Stanley's theorem does not
apply, e.g., the Tamari
lattices~\cite{ft:tfi,gra:lt,ht:spl} and
certain shuffle posets~\cite{gre:clm}.  Note also that we cannot drop
the level condition which replaced semimodularity completely:  If one
considers the non-crossing partition lattice then it has the same
modular chain as $\Pi_n$.  However, it does not satisfy the level
condition and its characteristic polynomial does not factor over
$\bbZ_{\ge0}$. 

I hope that you have enjoyed this tour through the world of the
characteristic polynomial and its factorizations.  Maybe you will feel
inspired to try one of the open problems mentioned along the way.

{\it Acknowledgment.}  I would like to thank the referee for very
helpful suggestions.

\begin{\bib}{99}

\bibitem{and:tp} G. Andrews, ``The Theory of Partitions,''
Addison-Wesley, Reading, MA, 1976.

\bibitem{ath:cps} C. A. Athanasiadis, Characteristic polynomials of
subspace arrangements and finite fields, {\it Adv.\ in Math.\ }
{\bf 122} (1966), 193--233.

\bibitem{al:sbr} C. A. Athanasiadis and S. Linusson, A simple
bijection for the regions of the Shi arrangement of hyperplanes, preprint.

\bibitem{bg:nbc}  H. Barcelo and A. Goupil, Non broken
circuits of reflection groups and their factorization in $D_n$,
{\it Israel J. Math.} {\bf 91} (1995), 285--306.

\bibitem{bs:gss} C. Bennett and B. E. Sagan, A generalization of
semimodular supersolvable lattices, {\it \jac}
{\bf 72} (1995), 209--231.

\bibitem{bjo:hsm} A. Bj\"orner, The homology and shellability of
matroids and geometric lattices, Chapter 7 in ``Matroid Applications,''  N.
White ed., Cambridge University Press, Cambridge, 1991, 226--283.

\bibitem{bjo:sa}  A. Bj\"orner, Subspace arrangements, in ``Proc. 1st
European Congress Math. (Paris 1992),'' A. Joseph and R. Rentschler
eds., Progress in Math., Vol.\ 122, Birkh\"auser, Boston, MA,
(1994), 321--370.

\bibitem{bl:ldt}  A. Bj\"orner and L. Lov\'asz,
Linear decision trees, subspace arrangements and M\"obius functions,
{\it \jams} {\bf 7} (1994), 667--706.

\bibitem{bly:ldt}  A. Bj\"orner, L. Lov\'asz and A. Yao,
Linear decision trees: volume estimates and topological bounds, in
``Proc. 24th ACM Symp. on Theory of Computing,''  ACM Press, New York,
NY, 1992, 170--177.

\bibitem{bs:sat}  A. Bj\"orner and B. Sagan, Subspace arrangements of
type $B_n$ and $D_n$, {\it \jac}, {\bf 5} (1996), 291--314.

\bibitem{bw:nsc1}  A. Bj\"orner and M. Wachs, Nonpure shellable
complexes and posets I, {\it Trans.\  Amer.\  Math.\  Soc.\ } {\bf 348}
(1996), 1299--1327.

\bibitem{bw:nsc2}  A. Bj\"orner and M. Wachs, Nonpure shellable 
complexes and posets II, {\it Trans.\ Amer.\ Math.\ Soc.\ }
{\bf 349} (1997), 3945--3975.

\bibitem{bw:hkm} A. Bj\"orner and V. Welker, The homology of
``$k$-equal'' manifolds and related partition lattices, {\it \aim}
{\bf 110} (1995), 277--313.

\bibitem{bs:mfl}  A. Blass and B. E. Sagan, M\"obius functions of
lattices, {\it Adv.\ in Math.} {\bf 127} (1997), 94--123.

\bibitem{bs:cep}  A. Blass and B. E. Sagan, Characteristic and
Ehrhart polynomials, {\it J. Algebraic Combin.} {\bf 7} (1998), 115--126.

\bibitem{bog:mfd} K. Bogart, The M\"obius function of the domination
lattice, unpublished manuscript, 1972.

\bibitem{bry:lip}  T. Brylawski, The lattice of integer partitions,
{\sl \dm} {\bf 6} (1973), 201--219.

\bibitem{bry:bcc}  T. Brylawski, The broken circuit complex, 
{\it \tams} {\bf 171} (1977), 235--282.

\bibitem{cl:gd} G. Chartrand and L. Lesniak, ``Graphs and Digraphs,'' 
second edition, Wadsworth \& Brooks/Cole, Pacific Grove, CA, 1986.

\bibitem{cr:cg}   H. Crapo and G.-C. Rota, ``On the Foundations of
Combinatorial Theory: Combinatorial Geometries,''  M.I.T. Press,
Cambridge, MA, 1970.

\bibitem{er:fha} P. H. Edelman and V. Reiner, Free hyperplane
arrangements between $A_{n-1}$ and $B_n$, {\it Math. Zeitschrift}
{\bf 215} (1994), 347--365.

\bibitem{ft:tfi} H. Friedman and D. Tamari, Probl\`emes
d'associativit\'e:  Une treillis finis induite par une loi
demi-associative, {\it J. Combin. Theory} {\bf 2} (1967), 215--242.

\bibitem{gra:lt} G. Gr\"atzer, ``Lattice Theory,'' Freeman and Co.,
San Francisco, CA, 1971, pp. 17--18, problems 26-36.

\bibitem{gre:clm} C. Greene, A class of lattices with M\"obius
function $\pm1$, {\sl \ejc} {\bf 9} (1988), 225--240.

\bibitem{gre:ps} C. Greene, Posets of Shuffles, {\sl \jcta} {\bf 47}
(1988), 191--206.

\bibitem{hai:cqr} M. Haiman, Conjectures on the quotient ring of
diagonal invariants. {\it J. Alg.\ Combin.\ }, {\bf 3} (1994),
17--76.

\bibitem{hea:rei} P. Headly, ``Reduced Expressions in Infinite Coxeter
Groups,'' Ph.D. thesis, University of Michigan, Ann Arbor, 1994.

\bibitem{hea:rwa} P. Headly, On reduced words in affine Weyl groups,
in ``Formal Power Series and Algebraic Combinatorics, May 23--27, 1994,'' 
DIMACS, Rutgers, 1994, 225--242.

\bibitem{ht:spl} S. Huang and D. Tamari, Problems of associativity:  A
simple proof for the lattice property of systems ordered by a
semi-associative law, {\it \jcta} {\bf 13} (1972), 7--13.

\bibitem{hum:rgc} J. E. Humphreys, ``Reflection Groups and Coxeter
Groups,'' Cambridge Studies in Advanced Mathematics, Cambridge
University Press, Cambridge, 1990.

\bibitem{jp:ifa} M. Jambu and L. Paris, Combinatorics of inductively
factored arrangements, {\it European J. Combin.} {\bf 16} (1995), 267--292.

\bibitem{js:bds} T. J\'ozefiak and B. E. Sagan, Basic derivations 
for subarrangements of Coxeter arrangements, {\it \jac} 
{\bf 2} (1993), 291--320.

\bibitem{kre:pnc}  G. Kreweras, Sur les partitions non-crois\'ees d'un
cycle, {\sl \dm} {\bf 1} (1972), 333--350.

\bibitem{lin:mfc}  S. Linusson, M\"obius functions and characteristic
polynomials for subspace arrangements embedded in $B_n$, preprint.

\bibitem{lin:prb}  S. Linusson, Partitions with restricted block
sizes, M\"obius functions and the $k$-of-each problem, 
{\it SIAM J. Discrete Math.\/} {\bf 10} (1997), 18--29.

\bibitem{ot:ah} P. Orlik and H. Terao, ``Arrangements of Hyperplanes,''
Grundlehren 300, Springer-Verlag, New York, NY, 1992.

\bibitem{pos:eag} A. Postnikov, ``Enumeration in algebra and
geometry,'' Ph.D. thesis, M.I.T., Cambridge, 1997. 

\bibitem{ps:dch} A. Postnikov and R. P. Stanley, Deformations of
Coxeter hyperplane arrangements, preprint.

\bibitem{rei:npc}  V. Reiner, Non-crossing partitions for classical
reflection groups, {\it Discrete Math.} {\bf 177} (1997), 195--222.

\bibitem{rot:tmf} G.-C. Rota, On the foundations of combinatorial
theory I. Theory of M\"obius functions, {\it Z. 
Wahrscheinlichkeitstheorie} {\bf 2} (1964), 340--368.

\bibitem{sai:tld} K. Saito, Theory of logarithmic differential
forms and logarithmic vector fields, {\it J. Fac. Sci. Univ. Tokyo
Sec. 1A Math.} {\bf 27} (1980), 265--291.

\bibitem{shi:klc} J. Y. Shi, The Kazhdan-Lusztig cells in certain
affine Weyl groups, Lecture Notes in Math., Vol. 1179, Springer-Verlag,
New York, NY, 1986.

\bibitem{shi:stc} J. Y. Shi, Sign types corresponding to an affine
Weyl group, {\it J. London Math.\ Soc.\ } {\bf 35} (1987), 56--74.

\bibitem{sta:sl}  R. P. Stanley, Supersolvable lattices, {\it
Algebra Universalis} {\bf 2} (1972), 197--217.

\bibitem{sta:ec1} R. P. Stanley, ``Enumerative Combinatorics,
Volume 1,''  Cambridge University Press, Cambridge, 1997.

\bibitem{sta:hai} R. P. Stanley, Hyperplane arrangements, interval
orders, and trees, {\it Proc.\ Nat.\ Acad.\ Sci.\/} {\bf 93} (1996)
2620--2625.

\bibitem{sun:aht} S. Sundaram, Applications of the Hopf trace formula
to computing homology representations, {\it Contemp.\ Math.} {\bf 178}
(1994), 277--309.

\bibitem{sw:hrk} S. Sundaram and M. Wachs, The homology
representations of the $k$-equal partition lattice, {\it \tams}
{\bf 349} (1997) 935--954.

\bibitem{sw:gal} S. Sundaram and V. Welker, Group actions on linear
subspace arrangements and applications to configuration spaces,
{\it \tams} {\bf 349} (1997) 1389--1420.

\bibitem{ter:ahf} H. Terao, Arrangements of hyperplanes and their
freeness I, II, {\it J. Fac. Sci. Univ. Tokyo}, {\bf 27} (1980), 293--320.

\bibitem{ter:gef} H. Terao, Generalized exponents of a free
arrangement of hyperplanes and the Shepherd-Todd-Brieskorn formula,
{\it \im} {\bf 63} (1981), 159--179.

\bibitem{ter:fah} H. Terao, Free arrangements of hyperplanes over an
arbitrary field, {\it \pja} {\bf 59} (1983), 301--303.

\bibitem{ter:jdf} H. Terao, The Jacobians and the discriminants of
finite reflection groups, {\it T\^ohoku Math. J.} {\bf 41} (1989),
237--247. 

\bibitem{ter:fos} H. Terao, Factorizations of Orlik-Solomon algebras,
{\it \aim} {\bf 91} (1992), 45--53.

\bibitem{wil:gen} H. S. Wilf, ``Generatingfunctionology,'' Academic
Press, Boston, MA, 1990.

\bibitem{zas:fua} T. Zaslavsky, ``Facing up to arrangements: Face-count
formulas for partitions of space by hyperplanes,''
Memoirs Amer. Math. Soc., No. 154, Amer. Math. Soc., Providence, RI, 1975.

\bibitem{zas:grs} T. Zaslavsky, The geometry of root systems
and signed graphs, {\it Amer. Math. Monthly} {\bf 88} (1981), 88--105.

\bibitem{zas:sgc} T. Zaslavsky, Signed graph coloring, {\it \dm}  
{\bf 39} (1982) 215--228.

\bibitem{zas:cis} T. Zaslavsky, Chromatic invariants of signed graphs,
{\it \dm}  {\bf 42} (1982) 287--312.

\bibitem{zie:ach} G. Ziegler, Algebraic combinatorics of hyperplane
arrangements, Ph.\ D. thesis, M.I.T., Cambridge, MA, 1987.

\end{\bib}

\end{document}